\documentclass[12 pt]{amsart}
\usepackage{amsfonts}
\usepackage{ifthen}
\usepackage{amsthm}
\usepackage{amsmath}
\usepackage{graphicx}
\usepackage{amscd,amssymb,amsthm}
\usepackage{graphicx}
\usepackage{epstopdf}
\usepackage{hyperref}

\newcounter{minutes}
\divide\time by 60
\newcounter{hours}
\multiply\time by 60 \addtocounter{minutes}{-\time}

\newtheorem{lemma}{Lemma}
\newtheorem{theorem}{Theorem}

\keywords{analytic function; $q$-Bessel functions; univalent functions; partial sums.}
\subjclass[2010]{30C45}

\title{On partial sums of normalized $q$-Bessel functions}

\author[H. Orhan]{Hal\.{i}t Orhan}
\address{Department of Mathematics, Faculty of Science, Atat\"{u}rk University, Erzurum, Turkey}
\email{orhanhalit607@gmail.com}
\author[\.{I}. Akta\c{s}]{\.{I}brah\.{I}m Akta\c{s}}
\address{Department of Mathematical Engineering, Faculty of Engineering and Natural Sciences, G\"{u}m\"{u}\c{s}hane University, G\"{u}m\"{u}\c{s}hane, Turkey}
\email{aktasibrahim38@gmail.com}

\begin{document}
	
	\def\thefootnote{}
	\footnotetext{ \texttt{File:~\jobname .tex,
			printed: \number\year-\number\month-\number\day,
			\thehours.\ifnum\theminutes<10{0}\fi\theminutes}
	} \makeatletter\def\thefootnote{\@arabic\c@footnote}\makeatother
	
	\maketitle
	
\begin{abstract}
	In the present investigation our main aim is to give lower bounds for the ratio of some normalized $q$-Bessel functions and their sequences of partial sums. Especially, we consider Jackson's second and third $q$-Bessel functions and we apply one normalization for each of them. 		
\end{abstract}

\section{Introduction}
Let $\mathcal{A}$ denote the class of functions of the following form: 
\begin{equation}\label{eq1.1}
f(z)=z+\sum_{n=2}^{\infty}a_{n}z^{n},
\end{equation}
which are analytic in the open unit disk  $$\mathcal{U}=\{z:z\in\mathbb{C}\text{ and }\left|z\right|<1\}.$$
We denote by $\mathcal{S}$ the class of all functions in  $\mathcal{A}$ which are univalent in $\mathcal{U}$.

The Jackson's second and third $q$-Bessel functions are defined by (see \cite{ANNABY})
\begin{equation}\label{eq1.2}
J_{\nu}^{(2)}(z;q)=\frac{(q^{\nu+1};q)_{\infty}}{(q;q)_{\infty}}\sum_{n\geq0}\frac{(-1)^{n}\left(\frac{z}{2}\right)^{2n+\nu}}{(q;q)_{n}(q^{\nu+1};q)_{n}}q^{n(n+\nu)}
\end{equation}
and
\begin{equation}\label{eq1.3}
J_{\nu}^{(3)}(z;q)=\frac{(q^{\nu+1};q)_{\infty}}{(q;q)_{\infty}}\sum_{n\geq0}\frac{(-1)^{n}z^{2n+\nu}}{(q;q)_{n}(q^{\nu+1};q)_{n}}q^{\frac{1}{2}{n(n+1)}},
\end{equation} where $z\in\mathbb{C},\nu>-1,q\in(0,1)$ and $$(a;q)_0=1,	(a;q)_n=\prod_{k=1}^{n}\left(1-aq^{k-1}\right),	 (a,q)_{\infty}=\prod_{k\geq1}\left(1-aq^{k-1}\right).$$
Here we would like to say that Jackson's third $q$-Bessel function is also known as Hahn-Exton $q$-Bessel function. 

Recently, the some geometric properties like univalence, starlikeness and convexity of the some special functions were investigated by many authors. Especially, in \cite{AB,BDM,bdoy,OA} authors have studied on the starlikeness and convexity of the some normalized $q$-Bessel functions. In addition, the some lower bounds for the ratio of some special functions and their sequences of partial sums were given in \cite{AO2,CD,OY,RADU}. Morever, results related with partial sums of analytic functions can be found in \cite{AO1,OSS,Sheil,Silverman,Silvia} etc.

Motivated by the previous works on analytic and some special functions, in this paper our aim is to present some lower bounds for the ratio of normalized $q$-Bessel functions to their sequences of partial sums. 

Due to the functions defined by \eqref{eq1.2} and \eqref{eq1.3} do not belong to the class $\mathcal{A}$, we consider following normalized forms of the $q$-Bessel functions:
\begin{equation}\label{eq1.4}
h_{\nu}^{(2)}(z;q)=2^{\nu}c_{\nu}(q)z^{1-\frac{\nu}{2}}J_{\nu}^{(2)}(\sqrt{z};q)=\sum_{n\geq0}K_{n}z^{n+1}
\end{equation}
and
\begin{equation}\label{eq1.5}
h_{\nu}^{(3)}(z;q)=c_{\nu}(q)z^{1-\frac{\nu}{2}}J_{\nu}^{(3)}(\sqrt{z};q)=\sum_{n\geq0}T_{n}z^{n+1},
\end{equation}
where $K_{n}=\frac{(-1)^nq^{n(n+\nu)}}{4^n(q;q)_{n}(q^{\nu+1};q)_{n}}, T_{n}=\frac{(-1)^nq^{\frac{1}{2}n(n+1)}}{(q;q)_{n}(q^{\nu+1};q)_{n}}\text{ and } c_{\nu}(q)=(q;q)_{\infty}\big/(q^{\nu+1};q)_{\infty}$. As a result of the above normalizations, all of the above functions belong to the class $\mathcal{A}$.

\section{Main Results}
\setcounter{equation}{0}
The following lemmas will be required in order to derive our main results.

\begin{lemma}
	Let $q\in\left(0,1\right), \nu>-1 \text{ and } 4(1-q)(1-q^{\nu})>q^{\nu}$. Then the function $h_{\nu}^{(2)}(z;q)$ satisfies tne next two inequalities for $z\in\mathcal{U}$:
	\begin{equation}\label{eq2.1}
	\left|h_{\nu}^{(2)}(z;q)\right|\leq\frac{4(1-q)(1-q^{\nu})}{4(1-q)(1-q^{\nu})-q^{\nu}},
	\end{equation}
	\begin{equation}\label{eq2.2}
	\left|\left(h_{\nu}^{(2)}(z;q)\right)^{\prime}\right|\leq\left(\frac{4(1-q)(1-q^{\nu})}{4(1-q)(1-q^{\nu})-q^{\nu}}\right)^{2}.
	\end{equation}
\end{lemma}

\begin{proof}
	It can be easily shown that the inequalities $$q^{n(n+\nu)}\leq q^{n\nu}, (1-q)^n\leq(q;q)_{n}\text{ and } (1-q^\nu)^n\leq(q^{\nu+1};q)_{n}$$ are valid for $q\in\left(0,1\right) \text{ and } \nu>-1$. Making use the above inequalities and well-known triangle inequality, for $z\in\mathcal{U}$, we get
	\begin{align*}
	\left|h_{\nu}^{(2)}(z;q)\right|&=\left|z+\sum_{n\geq1}\frac{(-1)^nq^{n(n+\nu)}}{4^n(q;q)_{n}(q^{\nu+1};q)_{n}}z^{n+1}\right|\\&\leq1+\sum_{n\geq1}\frac{q^{n(n+\nu)}}{4^n(q;q)_{n}(q^{\nu+1};q)_{n}}\\&\leq1+\sum_{n\geq1}\left(\frac{q^\nu}{4(1-q)(1-q^{\nu})}\right)^n\\&=1+\frac{q^\nu}{4(1-q)(1-q^{\nu})}\sum_{n\geq1}\left(\frac{q^\nu}{4(1-q)(1-q^{\nu})}\right)^{n-1}\\&=\frac{4(1-q)(1-q^{\nu})}{4(1-q)(1-q^{\nu})-q^{\nu}}
	\end{align*}
	and
	\begin{align*}
	\left|\left(h_{\nu}^{(2)}(z;q)\right)^{\prime}\right|&=\left|1+\sum_{n\geq1}\frac{(-1)^n(n+1)q^{n(n+\nu)}}{4^n(q;q)_{n}(q^{\nu+1};q)_{n}}z^{n}\right|\\&\leq1+\sum_{n\geq1}\frac{(n+1)q^{n(n+\nu)}}{4^n(q;q)_{n}(q^{\nu+1};q)_{n}}\\&\leq1+\sum_{n\geq1}(n+1)\left(\frac{q^\nu}{4(1-q)(1-q^{\nu})}\right)^n\\&=\left(\frac{4(1-q)(1-q^{\nu})}{4(1-q)(1-q^{\nu})-q^{\nu}}\right)^{2}.
	\end{align*}
	Thus, the inequalities \eqref{eq2.1} and \eqref{eq2.2} are proved.
\end{proof}

\begin{lemma}
	Let $q\in\left(0,1\right), \nu>-1 \text{ and } (1-q)(1-q^{\nu})>\sqrt{q}$. Then the function $h_{\nu}^{(3)}(z;q)$ satisfies the inequalities 
	\begin{equation}\label{eq2.3}
	\left|h_{\nu}^{(3)}(z;q)\right|\leq\frac{(1-q)(1-q^{\nu})}{(1-q)(1-q^{\nu})-\sqrt{q}},
	\end{equation}
	and
	\begin{equation}\label{eq2.4}
	\left|\left(h_{\nu}^{(3)}(z;q)\right)^{\prime}\right|\leq\left(\frac{(1-q)(1-q^{\nu})}{(1-q)(1-q^{\nu})-\sqrt{q}}\right)^{2}
	\end{equation}
	for $z\in\mathcal{U}.$
\end{lemma}

\begin{proof}
	It is known that the inequalities
	$$q^{\frac{1}{2}n(n+1)}\leq q^{\frac{1}{2}n}, (1-q)^n\leq(q;q)_{n}\text{ and } (1-q^\nu)^n\leq(q^{\nu+1};q)_{n}$$
	are valid for $q\in\left(0,1\right) \text{ and } \nu>-1$. Now, using the well-known triangle inequality for $z\in\mathcal{U}$, we have
	\begin{align*}
	\left|h_{\nu}^{(3)}(z;q)\right|&=\left|z+\sum_{n\geq1}\frac{(-1)^nq^{\frac{1}{2}n(n+1)}}{(q;q)_{n}(q^{\nu+1};q)_{n}}z^{n+1}\right|\\&\leq1+\sum_{n\geq1}\frac{q^{\frac{1}{2}n}}{(1-q)^n(1-q^\nu)^{n}}\\&\leq1+\frac{{\sqrt{q}}}{(1-q)(1-q^\nu)}\sum_{n\geq1}\left(\frac{{\sqrt{q}}}{(1-q)(1-q^\nu)}\right)^{n-1}\\&=\frac{(1-q)(1-q^{\nu})}{(1-q)(1-q^{\nu})-\sqrt{q}}
	\end{align*}
	and
	\begin{align*}
	\left|\left(h_{\nu}^{(3)}(z;q)\right)^{\prime}\right|&=\left|1+\sum_{n\geq1}\frac{(-1)^n(n+1)q^{\frac{1}{2}n(n+1)}}{(q;q)_{n}(q^{\nu+1};q)_{n}}z^{n}\right|\\&\leq1+\sum_{n\geq1}(n+1)\frac{q^{\frac{1}{2}n}}{(1-q)^n(1-q^\nu)^{n}}\\&\leq1+\frac{{\sqrt{q}}}{(1-q)(1-q^\nu)}\sum_{n\geq1}(n+1)\left(\frac{{\sqrt{q}}}{(1-q)(1-q^\nu)}\right)^{n-1}\\&=\left(\frac{(1-q)(1-q^{\nu})}{(1-q)(1-q^{\nu})-\sqrt{q}}\right)^2.
	\end{align*}
	So, the inequalities \eqref{eq2.3} and \eqref{eq2.4} are proved.
\end{proof}

Let $w(z)$ denote an analytic function in $\mathcal{U}$. In the proof of our main results, the following well-known result will be used frequently:
$$\Re\Bigg\{\frac{1+w(z)}{1-w(z)}\Bigg\}>0, \text{ if and only if }\left|w(z)\right|<1, z\in\mathcal{U}.$$

\begin{theorem}\label{th1}
	Let $\nu>-1, q\in\left(0,1\right),$ the function $h_{\nu}^{(2)}:\mathcal{U}\rightarrow\mathbb{C}$ be defined by \eqref{eq1.4} and its sequences of partial sums by $(h_{\nu}^{(2)})_{m}(z;q)=z+\sum_{n=1}^{m}K_{n}z^{n+1}.$ If the inequality $2(1-q)(1-q^{\nu})\geq{q^{\nu}},$ then the following inequalities hold true for $z\in\mathcal{U}$:
	\begin{equation}\label{eq2.5}
	\Re\Bigg\{\frac{h_{\nu}^{(2)}(z;q)}{(h_{\nu}^{(2)})_{m}(z;q)}\Bigg\}\geq\frac{4(1-q)(1-q^{\nu})-2q^{\nu}}{4(1-q)(1-q^{\nu})-q^{\nu}},
	\end{equation}
	\begin{equation}\label{eq2.6}
	\Re\Bigg\{\frac{(h_{\nu}^{(2)})_{m}(z;q)}{h_{\nu}^{(2)}(z;q)}\Bigg\}\geq\frac{4(1-q)(1-q^{\nu})-q^{\nu}}{4(1-q)(1-q^{\nu})}.
	\end{equation}
\end{theorem} 

\begin{proof}
	From the inequality \eqref{eq2.1} we have that
	\begin{equation}\label{eq2.7}
	1+\sum_{n\geq1}\left|K_{n}\right|\leq\frac{4(1-q)(1-q^{\nu})}{4(1-q)(1-q^{\nu})-q^{\nu}}.
	\end{equation}
    The inequalitiy \eqref{eq2.7} is equivalent to
	\begin{equation}\label{eq2.8}
	\frac{4(1-q)(1-q^{\nu})-q^{\nu}}{q^{\nu}}\sum_{n\geq1}\left|K_{n}\right|\leq1.
	\end{equation}
	In order to prove the inequality \eqref{eq2.5}, we consider the function $w(z)$ defined by $$\frac{1+w(z)}{1-w(z)}=\frac{4(1-q)(1-q^{\nu})-q^{\nu}}{q^{\nu}}\Bigg\{\frac{h_{\nu}^{(2)}(z;q)}{(h_{\nu}^{(2)})_{m}(z;q)}-\frac{4(1-q)(1-q^{\nu})-2q^{\nu}}{4(1-q)(1-q^{\nu})-q^{\nu}}\Bigg\}$$
	which is equivalent to 
	\begin{equation}\label{eq2.9}
	\frac{1+w(z)}{1-w(z)}=\frac{1+\sum_{n=1}^{m}K_{n}z^{n}+\frac{4(1-q)(1-q^{\nu})-q^{\nu}}{q^{\nu}}\sum_{n=m+1}^{\infty}K_{n}z^{n}}{1+\sum_{n=1}^{m}K_{n}z^{n}}.
	\end{equation}
	By using the equality \eqref{eq2.9} we get $$w(z)=\frac{\frac{4(1-q)(1-q^{\nu})-q^{\nu}}{q^{\nu}}\sum_{n=m+1}^{\infty}K_{n}z^{n}}{2+2\sum_{n=1}^{m}K_{n}z^{n}+\frac{4(1-q)(1-q^{\nu})-q^{\nu}}{q^{\nu}}\sum_{n=m+1}^{\infty}K_{n}z^{n}}$$ and $$\left|w(z)\right|\leq\frac{\frac{4(1-q)(1-q^{\nu})-q^{\nu}}{q^{\nu}}\sum_{n=m+1}^{\infty}\left|K_{n}\right|}{2-2\sum_{n=1}^{m}\left|K_{n}\right|-\frac{4(1-q)(1-q^{\nu})-q^{\nu}}{q^{\nu}}\sum_{n=m+1}^{\infty}\left|K_{n}\right|}.$$
	The inequality
	\begin{equation}\label{eq2.10}
	\sum_{n=1}^{m}\left|K_{n}\right|+\frac{4(1-q)(1-q^{\nu})-q^{\nu}}{q^{\nu}}\sum_{n=m+1}^{\infty}\left|K_{n}\right|\leq1
	\end{equation}
	implies that $\left|w(z)\right|\leq1$. It suffices to show that the left hand side of \eqref{eq2.10} is bounded above by $$\frac{4(1-q)(1-q^{\nu})-q^{\nu}}{q^{\nu}}\sum_{n\geq1}\left|K_{n}\right|,$$ which is equivalent to $$\frac{4(1-q)(1-q^{\nu})-2q^{\nu}}{q^{\nu}}\sum_{n\geq1}\left|K_{n}\right|\geq0.$$ The last inequality holds true for $2(1-q)(1-q^{\nu})\geq{q^{\nu}}.$
	
	In order to prove the result \eqref{eq2.6} we use the same method. Now, consider the function $p(z)$ given by
	$$\frac{1+p(z)}{1-p(z)}=\left(1+\frac{4(1-q)(1-q^{\nu})-q^{\nu}}{q^{\nu}}\right)\Bigg\{\frac{(h_{\nu}^{(2)})_{m}(z;q)}{h_{\nu}^{(2)}(z;q)}-\frac{4(1-q)(1-q^{\nu})-q^{\nu}}{4(1-q)(1-q^{\nu})}\Bigg\}.$$ Then from the last equality we get
	$$p(z)=\frac{-\frac{4(1-q)(1-q^{\nu})}{q^{\nu}}\sum_{n=m+1}^{\infty}K_{n}z^{n}}{2+2\sum_{n=1}^{m}K_{n}z^{n}-\frac{4(1-q)(1-q^{\nu})}{q^{\nu}}\sum_{n=m+1}^{\infty}K_{n}z^{n}}$$
	and
	$$\left|p(z)\right|\leq\frac{\frac{4(1-q)(1-q^{\nu})}{q^{\nu}}\sum_{n=m+1}^{\infty}\left|K_{n}\right|}{2-2\sum_{n=1}^{m}\left|K_{n}\right|-\frac{4(1-q)(1-q^{\nu})}{q^{\nu}}\sum_{n=m+1}^{\infty}\left|K_{n}\right|}.$$
	The inequality
	\begin{equation}\label{eq2.11}
	\sum_{n=1}^{m}\left|K_{n}\right|+\frac{4(1-q)(1-q^{\nu})}{q^{\nu}}\sum_{n=m+1}^{\infty}\left|K_{n}\right|\leq1
	\end{equation}
	implies that $\left|p(z)\right|\leq1$. Since the left hand side of \eqref{eq2.11} is bounded above by $$\frac{4(1-q)(1-q^{\nu})-q^{\nu}}{q^{\nu}}\sum_{n=1}^{m}\left|K_{n}\right|\geq0$$ the proof is completed.
\end{proof}

\begin{theorem}\label{th2}
	Let $\nu>-1, q\in\left(0,1\right),$ the function $h_{\nu}^{(2)}:\mathcal{U}\rightarrow\mathbb{C}$ be defined by \eqref{eq1.4} and its sequences of partial sums by $(h_{\nu}^{(2)})_{m}(z;q)=z+\sum_{n=1}^{m}K_{n}z^{n+1}.$ If the inequality $ (1-q)(1-q^{\nu})\geq{q^{\nu}} $ is valid, then the following inequalities hold true for $z\in\mathcal{U}$:
	\begin{equation}\label{eq2.12}
	\Re\Bigg\{\frac{\left(h_{\nu}^{(2)}(z;q)\right)^{\prime}}{\left((h_{\nu}^{(2)})_{m}(z;q)\right)^{\prime}}\Bigg\}\geq\frac{16(1-q)(1-q^{\nu})\left((1-q)(1-q^{\nu})-q^{\nu}\right)+2q^{2\nu}}{8(1-q)(1-q^{\nu})q^{\nu}-q^{2\nu}},
	\end{equation}
	\begin{equation}\label{eq2.13}
	\Re\Bigg\{\frac{\left((h_{\nu}^{(2)})_{m}(z;q)\right)^{\prime}}{\left(h_{\nu}^{(2)}(z;q)\right)^{\prime}}\Bigg\}\geq\frac{\left(4(1-q)(1-q^{\nu})-q^{\nu}\right)^{2}}{8(1-q)(1-q^{\nu})q^{\nu}-q^{2\nu}}.
	\end{equation}
\end{theorem}
\begin{proof}
	From the inequality \eqref{eq2.2} we have that
	\begin{equation}\label{eq2.14}
	1+\sum_{n\geq1}\left(n+1\right)\left|K_{n}\right|\leq\left(\frac{4(1-q)(1-q^{\nu})}{4(1-q)(1-q^{\nu})-q^{\nu}}\right)^2.
	\end{equation}
	The inequality \eqref{eq2.14} is equivalent to
	\begin{equation}\label{eq2.15}
	\frac{\left(4(1-q)(1-q^{\nu})-q^{\nu}\right)^2}{8(1-q)(1-q^\nu)q^{\nu}-q^{2\nu}}\sum_{n\geq1}\left(n+1\right)\left|K_{n}\right|\leq1.
	\end{equation}
	In order to prove the inequality \eqref{eq2.12}, we consider the function $h(z)$ defined by $$\frac{1+h(z)}{1-h(z)}=\frac{\left(4(1-q)(1-q^{\nu})-q^{\nu}\right)^2}{8(1-q)(1-q^\nu)q^{\nu}-q^{2\nu}}\Bigg\{\frac{\left(h_{\nu}^{(2)}(z;q)\right)^{\prime}}{\left((h_{\nu}^{(2)})_{m}(z;q)\right)^{\prime}} -\delta\Bigg\},$$ where $\delta=\frac{16(1-q)(1-q^{\nu})\left((1-q)(1-q^{\nu})-q^{\nu}\right)+2q^{2\nu}}{8(1-q)(1-q^{\nu})q^{\nu}-q^{2\nu}}.$
	The last equality is equivalent to	
	\begin{equation}\label{eq2.16}
	\frac{1+h(z)}{1-h(z)}=\frac{1+{\sum_{n=1}^{m}}\left(n+1\right)K_{n}z^{n}+\frac{\left(4(1-q)(1-q^{\nu})-q^{\nu}\right)^2}{8(1-q)(1-q^\nu)q^{\nu}-q^{2\nu}}\sum_{n=m+1}^{\infty}\left(n+1\right)K_{n}z^{n}}{1+\sum_{n=1}^{m}\left(n+1\right)K_{n}z^{n}}.
	\end{equation}
	By using the equality \eqref{eq2.16} we get 
	$$h(z)=\frac{\frac{\left(4(1-q)(1-q^{\nu})-q^{\nu}\right)^2}{8(1-q)(1-q^\nu)q^{\nu}-q^{2\nu}}\sum_{n=m+1}^{\infty}\left(n+1\right)K_{n}z^{n}}{2+2\sum_{n=1}^{m}\left(n+1\right)K_{n}z^{n}+\frac{\left(4(1-q)(1-q^{\nu})-q^{\nu}\right)^2}{8(1-q)(1-q^\nu)q^{\nu}-q^{2\nu}}\sum_{n=m+1}^{\infty}\left(n+1\right)K_{n}z^{n}}$$ 
	and 
	$$\left|h(z)\right|\leq\frac{\frac{\left(4(1-q)(1-q^{\nu})-q^{\nu}\right)^2}{8(1-q)(1-q^\nu)q^{\nu}-q^{2\nu}}\sum_{n=m+1}^{\infty}\left(n+1\right)\left|K_{n}\right|}{2-2\sum_{n=1}^{m}\left(n+1\right)\left|K_{n}\right|-\frac{\left(4(1-q)(1-q^{\nu})-q^{\nu}\right)^2}{8(1-q)(1-q^\nu)q^{\nu}-q^{2\nu}}\sum_{n=m+1}^{\infty}\left(n+1\right)\left|K_{n}\right|}.$$
	The inequality
	\begin{equation}\label{eq2.17}
	\sum_{n=1}^{m}\left(n+1\right)\left|K_{n}\right|+\frac{\left(4(1-q)(1-q^{\nu})-q^{\nu}\right)^2}{8(1-q)(1-q^\nu)q^{\nu}-q^{2\nu}}\sum_{n=m+1}^{\infty}\left(n+1\right)\left|K_{n}\right|\leq1
	\end{equation}
	implies that $\left|h(z)\right|\leq1$. It suffices to show that the left hand side of \eqref{eq2.17} is bounded above by $$\frac{\left(4(1-q)(1-q^{\nu})-q^{\nu}\right)^2}{8(1-q)(1-q^\nu)q^{\nu}-q^{2\nu}}\sum_{n\geq1}\left(n+1\right)\left|K_{n}\right|,$$ which is equivalent to $$\delta\sum_{n\geq1}\left(n+1\right)\left|K_{n}\right|\geq0.$$ Thus, the result \eqref{eq2.12} is proved.
	
	To prove the result \eqref{eq2.13}, consider the function $k(z)$ defined by
	$$\frac{1+k(z)}{1-k(z)}=\Bigg\{1+\frac{\left(4(1-q)(1-q^{\nu})-q^{\nu}\right)^2}{8(1-q)(1-q^\nu)q^{\nu}-q^{2\nu}}\Bigg\}\Bigg\{\frac{\left(h_{\nu}^{(2)}(z;q)\right)^{\prime}}{\left((h_{\nu}^{(2)})_{m}(z;q)\right)^{\prime}} -\frac{\left(4(1-q)(1-q^{\nu})-q^{\nu}\right)^2}{8(1-q)(1-q^\nu)q^{\nu}-q^{2\nu}}\Bigg\}.$$ The last equality is equivalent to
	\begin{equation}\label{eq2.18}
	\frac{1+k(z)}{1-k(z)}=\frac{1+\sum_{n=1}^{m}(n+1)K_{n}z^{n}-\frac{\left(4(1-q)(1-q^{\nu})-q^{\nu}\right)^2}{8(1-q)(1-q^\nu)q^{\nu}-q^{2\nu}}\sum_{n=m+1}^{\infty}(n+1)K_{n}z^{n}}{1+\sum_{n\geq1}(n+1)K_{n}z^{n}}.
	\end{equation}
	From the equality \eqref{eq2.17} we have
	$$k(z)=\frac{-\frac{16(1-q)^2(1-q^{\nu})^2}{8(1-q)(1-q^\nu)q^{\nu}-q^{2\nu}}\sum_{n=m+1}^{\infty}(n+1)K_{n}z^{n}}{2+2\sum_{n=1}^{m}(n+1)K_{n}z^{n}-\delta\sum_{n=m+1}^{\infty}(n+1)K_{n}z^{n}}$$ and 
	$$\left|k(z)\right|\leq\frac{\frac{16(1-q)^2(1-q^{\nu})^2}{8(1-q)(1-q^\nu)q^{\nu}-q^{2\nu}}\sum_{n=m+1}^{\infty}(n+1)\left|K_{n}\right|}{2-2\sum_{n=1}^{m}(n+1)\left|K_{n}\right|-\delta\sum_{n=m+1}^{\infty}(n+1)\left|K_{n}\right|}.$$
	The inequality
	\begin{equation}\label{eq2.19}
	\sum_{n=1}^{m}(n+1)\left|K_{n}\right|+\frac{\left(4(1-q)(1-q^{\nu})-q^{\nu}\right)^2}{8(1-q)(1-q^\nu)q^{\nu}-q^{2\nu}}\sum_{n=m+1}^{\infty}(n+1)\left|K_{n}\right|\leq1
	\end{equation}
	implies that $\left|k(z)\right|\leq1.$ Since the left hand side of \eqref{eq2.19} is bounded above by
	$$\frac{\left(4(1-q)(1-q^{\nu})-q^{\nu}\right)^2}{8(1-q)(1-q^\nu)q^{\nu}-q^{2\nu}}\sum_{n\geq1}(n+1)\left|K_{n}\right|,$$ which is equivalent to 
	$$\delta\sum_{n=m+1}^{\infty}(n+1)\left|K_{n}\right|\geq0,$$
	the proof of result \eqref{eq2.13} is completed.
\end{proof}

\begin{theorem}\label{th3}
	Let $\nu>-1, q\in\left(0,1\right),$ the function $h_{\nu}^{(3)}:\mathcal{U}\rightarrow\mathbb{C}$ be defined by \eqref{eq1.5} and its sequences of partial sums by $(h_{\nu}^{(3)})_{m}(z;q)=z+\sum_{n=1}^{m}T_{n}z^{n+1}.$ If the inequality $(1-q)(1-q^{\nu})\geq2\sqrt{q}$ is valid, then the next two inequalities are valid for $z\in\mathcal{U}$:
	\begin{equation}\label{eq2.20}
	\Re\Bigg\{\frac{h_{\nu}^{(3)}(z;q)}{(h_{\nu}^{(3)})_{m}(z;q)}\Bigg\}\geq\frac{(1-q)(1-q^{\nu})-2\sqrt{q}}{\sqrt{q}},
	\end{equation}
	\begin{equation}\label{eq2.21}
	\Re\Bigg\{\frac{(h_{\nu}^{(3)})_{m}(z;q)}{h_{\nu}^{(3)}(z;q)}\Bigg\}\geq\frac{(1-q)(1-q^{\nu})-\sqrt{q}}{\sqrt{q}}.
	\end{equation}
\end{theorem}

\begin{proof}
	From the inequality \eqref{eq2.3} we have that
	\begin{equation}\label{eq2.22}
	1+\sum_{n\geq1}\left|T_{n}\right|\leq\frac{(1-q)(1-q^{\nu})}{(1-q)(1-q^{\nu})-\sqrt{q}}.
	\end{equation}
	The inequalitiy \eqref{eq2.22} is equivalent to
	\begin{equation}\label{eq2.23}
	\frac{(1-q)(1-q^{\nu})-\sqrt{q}}{\sqrt{q}}\sum_{n\geq1}\left|T_{n}\right|\leq1.
	\end{equation}
	In order to prove the inequality \eqref{eq2.20}, we consider the function $\phi(z)$ defined by $$\frac{1+\phi(z)}{1-\phi(z)}=\frac{(1-q)(1-q^{\nu})-\sqrt{q}}{\sqrt{q}}\Bigg\{\frac{h_{\nu}^{(3)}(z;q)}{(h_{\nu}^{(3)})_{m}(z;q)}-\frac{(1-q)(1-q^{\nu})-2\sqrt{q}}{\sqrt{q}}\Bigg\},$$
	which is equivalent to 
	\begin{equation}\label{eq2.24}
	\frac{1+\phi(z)}{1-\phi(z)}=\frac{1+\sum_{n=1}^{m}T_{n}z^{n}+\frac{(1-q)(1-q^{\nu})-\sqrt{q}}{\sqrt{q}}\sum_{n=m+1}^{\infty}T_{n}z^{n}}{1+\sum_{n=1}^{m}T_{n}z^{n}}.
	\end{equation}
	From the equality \eqref{eq2.24} we obtain $$\phi(z)=\frac{\frac{(1-q)(1-q^{\nu})-\sqrt{q}}{\sqrt{q}}\sum_{n=m+1}^{\infty}T_{n}z^{n}}{2+2\sum_{n=1}^{m}T_{n}z^{n}+\frac{(1-q)(1-q^{\nu})-\sqrt{q}}{\sqrt{q}}\sum_{n=m+1}^{\infty}T_{n}z^{n}}$$ and $$\left|\phi(z)\right|\leq\frac{\frac{(1-q)(1-q^{\nu})-\sqrt{q}}{\sqrt{q}}\sum_{n=m+1}^{\infty}\left|T_{n}\right|}{2-2\sum_{n=1}^{m}\left|T_{n}\right|-\frac{(1-q)(1-q^{\nu})-\sqrt{q}}{\sqrt{q}}\sum_{n=m+1}^{\infty}\left|T_{n}\right|}.$$
	The inequality
	\begin{equation}\label{eq2.25}
	\sum_{n=1}^{m}\left|T_{n}\right|+\frac{(1-q)(1-q^{\nu})-\sqrt{q}}{\sqrt{q}}\sum_{n=m+1}^{\infty}\left|T_{n}\right|\leq1
	\end{equation}
	implies that $\left|\phi(z)\right|\leq1$. It suffices to show that the left hand side of \eqref{eq2.25} is bounded above by $$\frac{(1-q)(1-q^{\nu})-\sqrt{q}}{\sqrt{q}}\sum_{n\geq1}\left|T_{n}\right|,$$ which is equivalent to $$\frac{(1-q)(1-q^{\nu})-2\sqrt{q}}{\sqrt{q}}\sum_{n=1}^{m}\left|T_{n}\right|\geq0.$$ The last inequality holds true for $(1-q)(1-q^{\nu})\geq{2\sqrt{q}}.$
	
	In order to prove the result \eqref{eq2.21}, we consider the function $\varphi(z)$ given by
	$$\frac{1+\varphi(z)}{1-\varphi(z)}=\left(1+\frac{(1-q)(1-q^{\nu})-\sqrt{q}}{\sqrt{q}}\right)\Bigg\{\frac{(h_{\nu}^{(3)})_{m}(z;q)}{h_{\nu}^{(3)}(z;q)}-\frac{(1-q)(1-q^{\nu})-\sqrt{q}}{\sqrt{q}}\Bigg\}.$$ Then from the last equality we get
	$$\varphi(z)=\frac{-\frac{(1-q)(1-q^{\nu})-\sqrt{q}}{\sqrt{q}}\sum_{n=m+1}^{\infty}T_{n}z^{n}}{2+2\sum_{n=1}^{m}T_{n}z^{n}-\frac{(1-q)(1-q^{\nu})-\sqrt{q}}{\sqrt{q}}\sum_{n=m+1}^{\infty}T_{n}z^{n}}$$
	and
	$$\left|\varphi(z)\right|\leq\frac{\frac{(1-q)(1-q^{\nu})-\sqrt{q}}{\sqrt{q}}\sum_{n=m+1}^{\infty}\left|T_{n}\right|}{2-2\sum_{n=1}^{m}\left|T_{n}\right|-\frac{(1-q)(1-q^{\nu})-\sqrt{q}}{\sqrt{q}}\sum_{n=m+1}^{\infty}\left|T_{n}\right|}.$$
	The inequality
	\begin{equation}\label{eq2.26}
	\sum_{n=1}^{m}\left|T_{n}\right|+\frac{(1-q)(1-q^{\nu})-\sqrt{q}}{\sqrt{q}}\sum_{n=m+1}^{\infty}\left|T_{n}\right|\leq1
	\end{equation}
	implies that $\left|\varphi(z)\right|\leq1$. Since the left hand side of \eqref{eq2.26} is bounded above by 
	$$\frac{(1-q)(1-q^{\nu})-\sqrt{q}}{\sqrt{q}}\sum_{n\geq1}\left|T_{n}\right|,$$
	which is equivalent to
	$$\frac{(1-q)(1-q^{\nu})-2\sqrt{q}}{\sqrt{q}}\sum_{n=1}^{m}\left|T_{n}\right|\geq0.$$This completes the proof of the theorem.
\end{proof}

\begin{theorem}\label{th4}
	Let $\nu>-1, q\in\left(0,1\right),$ the function $h_{\nu}^{(3)}:\mathcal{U}\rightarrow\mathbb{C}$ be defined by \eqref{eq1.5} and its sequences of partial sums by $(h_{\nu}^{(3)})_{m}(z;q)=z+\sum_{n=1}^{m}T_{n}z^{n+1}.$ If the inequality $(1-q)(1-q^{\nu})\geq4\sqrt{q},$ then the next two inequalities are valid for $z\in\mathcal{U}$:
	\begin{equation}\label{eq2.27}
	\Re\Bigg\{\frac{\left(h_{\nu}^{(3)}(z;q)\right)^{\prime}}{\left((h_{\nu}^{(3)})_{m}(z;q)\right)^{\prime}}\Bigg\}\geq\frac{(1-q)^2(1-q^{\nu})^2-4(1-q)(1-q^{\nu})\sqrt{q}+2q}{2(1-q)(1-q^{\nu})\sqrt{q}-q},
	\end{equation}
	\begin{equation}\label{eq2.28}
	\Re\Bigg\{\frac{\left((h_{\nu}^{(3)})_{m}(z;q)\right)^{\prime}}{\left(h_{\nu}^{(3)}(z;q)\right)^{\prime}}\Bigg\}\geq\frac{\left((1-q)(1-q^{\nu})-\sqrt{q}\right)^2}{2(1-q)(1-q^{\nu})\sqrt{q}-q}.
    \end{equation}
    \end{theorem}
    
\begin{proof}
		From the inequality \eqref{eq2.4} we have that
		\begin{equation}\label{eq2.29}
		1+\sum_{n\geq1}\left(n+1\right)\left|T_{n}\right|\leq\left(\frac{(1-q)(1-q^{\nu})}{(1-q)(1-q^{\nu})-\sqrt{q}}\right)^2.
		\end{equation}
		The inequality \eqref{eq2.29} is equivalent to
		\begin{equation}\label{eq2.30}
		\frac{\left((1-q)(1-q^{\nu})-\sqrt{q}\right)^2}{2(1-q)(1-q^\nu)\sqrt{q}-q}\sum_{n\geq1}\left(n+1\right)\left|T_{n}\right|\leq1.
		\end{equation}
		In order to prove the inequality \eqref{eq2.27}, we consider the function $\psi(z)$ defined by $$\frac{1+\psi(z)}{1-\psi(z)}=\frac{\left((1-q)(1-q^{\nu})-\sqrt{q}\right)^2}{2(1-q)(1-q^\nu)\sqrt{q}-q}\Bigg\{\frac{\left(h_{\nu}^{(3)}(z;q)\right)^{\prime}}{\left((h_{\nu}^{(3)})_{m}(z;q)\right)^{\prime}} -\lambda\Bigg\},$$ where $\lambda=\frac{(1-q)^{2}(1-q^{\nu})^{2}-4(1-q)(1-q^{\nu})\sqrt{q}+2q}{2(1-q)(1-q^{\nu})\sqrt{q}-q}.$
		The last equality is equivalent to	
		\begin{equation}\label{eq2.31}
		\frac{1+\psi(z)}{1-\psi(z)}=\frac{1+{\sum_{n=1}^{m}}\left(n+1\right)T_{n}z^{n}+\frac{\left((1-q)(1-q^{\nu})-\sqrt{q}\right)^2}{2(1-q)(1-q^\nu)\sqrt{q}-q}\sum_{n=m+1}^{\infty}\left(n+1\right)T_{n}z^{n}}{1+\sum_{n=1}^{m}\left(n+1\right)T_{n}z^{n}}.
		\end{equation}
		By using the equality \eqref{eq2.31} we get 
		$$\psi(z)=\frac{\frac{\left((1-q)(1-q^{\nu})-\sqrt{q}\right)^2}{2(1-q)(1-q^\nu)\sqrt{q}-q}\sum_{n=m+1}^{\infty}\left(n+1\right)T_{n}z^{n}}{2+2\sum_{n=1}^{m}\left(n+1\right)T_{n}z^{n}+\frac{\left((1-q)(1-q^{\nu})-\sqrt{q}\right)^2}{2(1-q)(1-q^\nu)\sqrt{q}-q}\sum_{n=m+1}^{\infty}\left(n+1\right)T_{n}z^{n}}$$ 
		and 
		$$\left|\psi(z)\right|\leq\frac{\frac{\left((1-q)(1-q^{\nu})-\sqrt{q}\right)^2}{2(1-q)(1-q^\nu)\sqrt{q}-q}\sum_{n=m+1}^{\infty}\left(n+1\right)\left|T_{n}\right|}{2-2\sum_{n=1}^{m}\left(n+1\right)\left|T_{n}\right|-\frac{\left((1-q)(1-q^{\nu})-\sqrt{q}\right)^2}{2(1-q)(1-q^\nu)\sqrt{q}-q}\sum_{n=m+1}^{\infty}\left(n+1\right)\left|T_{n}\right|}.$$
		The inequality
		\begin{equation}\label{eq2.32}
		\sum_{n=1}^{m}\left(n+1\right)\left|T_{n}\right|+\frac{\left((1-q)(1-q^{\nu})-\sqrt{q}\right)^2}{2(1-q)(1-q^\nu)\sqrt{q}-q}\sum_{n=m+1}^{\infty}\left(n+1\right)\left|T_{n}\right|\leq1
		\end{equation}
		implies that $\left|\psi(z)\right|\leq1$. It suffices to show that the left hand side of \eqref{eq2.32} is bounded above by $$\frac{\left((1-q)(1-q^{\nu})-\sqrt{q}\right)^2}{2(1-q)(1-q^\nu)\sqrt{q}-q}\sum_{n\geq1}\left(n+1\right)\left|T_{n}\right|,$$ which is equivalent to $$\lambda\sum_{n=1}^{m}\left(n+1\right)\left|T_{n}\right|\geq0.$$ Thus, the result \eqref{eq2.27} is proved.
		
		To prove the result \eqref{eq2.28}, consider the function $\rho(z)$ defined by
		$$\frac{1+\rho(z)}{1-\rho(z)}=\Bigg\{1+\frac{\left((1-q)(1-q^{\nu})-\sqrt{q}\right)^2}{2(1-q)(1-q^\nu)\sqrt{q}-q}\Bigg\}\Bigg\{\frac{\left(h_{\nu}^{(3)}(z;q)\right)^{\prime}}{\left((h_{\nu}^{(3)})_{m}(z;q)\right)^{\prime}} -\frac{\left((1-q)(1-q^{\nu})-\sqrt{q}\right)^2}{2(1-q)(1-q^\nu)\sqrt{q}-q}\Bigg\}.$$ The last equality is equivalent to
		\begin{equation}\label{eq2.33}
		\frac{1+\rho(z)}{1-\rho(z)}=\frac{1+\sum_{n=1}^{m}(n+1)T_{n}z^{n}-\frac{\left((1-q)(1-q^{\nu})-\sqrt{q}\right)^2}{2(1-q)(1-q^\nu)\sqrt{q}-q}\sum_{n=m+1}^{\infty}(n+1)T_{n}z^{n}}{1+\sum_{n=1}^{\infty}(n+1)T_{n}z^{n}}.
		\end{equation}
		From the equality \eqref{eq2.33} we get
		$$\rho(z)=\frac{-\frac{\left((1-q)(1-q^{\nu})-\sqrt{q}\right)^2}{2(1-q)(1-q^\nu)\sqrt{q}-q}\sum_{n=m+1}^{\infty}(n+1)T_{n}z^{n}}{2+2\sum_{n=1}^{m}(n+1)T_{n}z^{n}-\frac{\left((1-q)(1-q^{\nu})-\sqrt{q}\right)^2}{2(1-q)(1-q^\nu)\sqrt{q}-q}\sum_{n=m+1}^{\infty}(n+1)T_{n}z^{n}}$$
		and 
		$$\left|\rho(z)\right|\leq\frac{\frac{\left((1-q)(1-q^{\nu})-\sqrt{q}\right)^2}{2(1-q)(1-q^\nu)\sqrt{q}-q}\sum_{n=m+1}^{\infty}(n+1)\left|T_{n}\right|}{2-2\sum_{n=1}^{m}(n+1)\left|T_{n}\right|-\frac{\left((1-q)(1-q^{\nu})-\sqrt{q}\right)^2}{2(1-q)(1-q^\nu)\sqrt{q}-q}\sum_{n=m+1}^{\infty}(n+1)\left|T_{n}\right|}.$$
		The inequality
		\begin{equation}\label{eq2.34}
		\sum_{n=1}^{m}(n+1)\left|T_{n}\right|+\frac{\left((1-q)(1-q^{\nu})-\sqrt{q}\right)^2}{2(1-q)(1-q^\nu)\sqrt{q}-q}\sum_{n=m+1}^{\infty}(n+1)\left|T_{n}\right|\leq1
		\end{equation}
		implies that $\left|\rho(z)\right|\leq1.$ Since the left hand side of \eqref{eq2.34} is bounded above by
		$$\frac{\left((1-q)(1-q^{\nu})-\sqrt{q}\right)^2}{2(1-q)(1-q^\nu)\sqrt{q}-q}\sum_{n\geq1}(n+1)\left|T_{n}\right|,$$ which is equivalent to 
		$$\frac{(1-q)(1-q^{\nu})\left((1-q)(1-q^\nu)-4\sqrt{q}\right)+2q}{2(1-q)(1-q^\nu)\sqrt{q}-q}\sum_{n=1}^{m}(n+1)\left|T_{n}\right|\geq0,$$
		the proof of result \eqref{eq2.28} is completed.
\end{proof}

\end{document}